\documentclass[aip,jmp,graphicx]{revtex4-1}
\usepackage{amsfonts,amsthm,amsmath,amssymb,amscd}
\usepackage{esint}
\usepackage{graphics}
\usepackage{indentfirst}
\usepackage{latexsym}
\usepackage[dvips]{epsfig}
\usepackage{mathrsfs}
\usepackage{enumitem}
\usepackage{cases}
\usepackage{hyperref}
\usepackage{verbatim}

\allowdisplaybreaks
\numberwithin{equation}{section}

\newcommand{\Ubar}{\bar{U}}

\newcommand{\rhobar}{\bar{\rho}}
\newcommand{\ubar}{\bar{u}}
\newcommand{\sbar}{\bar{s}}
\newcommand{\pbar}{\bar{p}}
\newcommand{\cbar}{\bar{c}}

\newtheorem{theorem}{Theorem}[section]

\newtheorem{lemma}[theorem]{Lemma}

\newtheorem{remark}[theorem]{Remark}

\begin{document}

\title{Nonlinear kink modes for nonisentropic  supersonic vortex sheets}

\author{Zihao Li}
\email[]{zihaol21210@gmail.com}

\author{Yang Xu}
\email[Corresponding author: ]{youngxu0518@163.com}

\affiliation{School of Mathematical Sciences, South China Normal University, Guangzhou, 510631, China}

\begin{abstract}
In this paper, planar nonisentropic vortex sheets of the two-dimensional Euler equations under the stability condition are studied via nonlinear geometric optics, focusing on the interaction induced by a small-amplitude, high-frequency incident wave impinging on a planar vortex sheet. It is found that there exist three \emph{asymmetric} critical incident angles at which resonance occurs, leading to a significant amplification of the reflected and refracted oscillatory waves. Moreover, by refining the asymptotic expansion, we deduce that the leading profiles of the highly oscillatory waves are described by initial-boundary value problems for Burgers transport equations, while the leading profile of the contact discontinuity front satisfies an initial
value problem for a Hamilton--Jacobi-type equation, which reduces to a forced linear transport equation on one resonant branch. These results extend the resonant amplification mechanism from the isentropic setting to the nonisentropic Euler flow and reveal the \emph{asymmetric} influence of entropy variations on the nonlinear evolution of compressible vortex sheets.
\end{abstract}
\maketitle 
\pagestyle{plain}

\section{Introduction}
We consider the interaction between high-frequency oscillatory waves and vortex sheets for the two-dimensional nonisentropic compressible Euler equations of gas dynamics. \begin{equation}
\left\{
\begin{aligned}
   & D_t p+\gamma p\,\nabla\cdot\mathbf{u} = 0,\\
   & \rho D_t\mathbf{u}+\nabla p = 0,\\
   & D_t s = 0,
\end{aligned}
\right.
\end{equation}
where $D_t:=\partial_t+\mathbf{u}\cdot\nabla$, $\mathbf{u}=(u_1,u_2)$ denotes the velocity field, and $p$, $\rho$, and $s$ denote the pressure, density, and entropy, respectively. The study of discontinuous solutions and their interaction with oscillatory waves for such hyperbolic systems is a fundamental issue in gas dynamics and in the theory of nonlinear wave propagation.

A natural starting point for addressing these interaction problems is the well-posedness and stability theory for the underlying nonlinear waves and characteristic discontinuities. For fundamental wave patterns, such as uniformly stable shocks and multidimensional rarefaction waves, several important results in this direction have been established in \cite{Majda.1983a,Majda.1983}.

For the Euler equations, there have been many interesting results on the stability and nonlinear behavior of vortex sheets. In two space dimensions, Miles~\cite{Miles.1957, Miles.1958} showed by mode analysis that vortex sheets for the compressible isentropic Euler equations are linearly stable when the Mach number is larger than $\sqrt{2}$ and unstable when it is smaller than $\sqrt{2}$. The corresponding rigorous theory was later developed by Coulombel--Secchi. They first obtained the linear stability theory for two-dimensional isentropic compressible vortex sheets under the strict stability condition~\cite{Coulo.S2004}, and then proved short-time existence and nonlinear stability of the free-boundary problem in~\cite{Coulo.S2008}.

In the nonisentropic case, the linear stability of piecewise constant vortex sheets for the two-dimensional Euler equations was also analyzed by Miles~\cite{Miles.1958}; see also Coulombel--Morando~\cite{Coulo.Mora.2004} for an algebraic approach. Morando--Trebeschi~\cite{Moran.T2008} established the weak linear stability of two-dimensional nonisentropic vortex sheets.  Building on this weak linear stability result, Morando--Trebeschi--Wang~\cite{Morando.2019} established the well-posedness of the effective linear problem and proved the short-time existence and nonlinear stability of nonisentropic vortex sheets. Another related class of nonisentropic contact discontinuities is entropy waves, whose nonlinear stability was recently studied by Wang--Zhang--Zhao~\cite{Wang.Z.Z2024}; unlike vortex sheets, they have continuous velocity but allow jumps in density and entropy.

Beyond the two-dimensional Euler setting, related stability problems have also been studied for higher-dimensional Euler flows. A stability analysis of vortex sheets in three-dimensional gas dynamics was conducted by Fejer and Miles \cite{Fejer.M1963}, revealing an instability that echoes the Kelvin--Helmholtz instability in incompressible fluids. In a related three-dimensional steady Euler framework, Wang--Yu~\cite{Wang.Y2013a,Wang.Y2015} studied the weak stability and structural nonlinear stability of supersonic contact discontinuities.

Analogous interface stability problems have also been studied in magnetohydrodynamics (MHD). For magnetic fields tangent to the interface, the corresponding current-vortex sheets have been investigated in both compressible and incompressible flows (see~\cite{Moran.S.T.Y2023,Trakh.2009,Wang.Y2013,Moran.T.T2008}). For MHD contact discontinuities with a continuous magnetic field and a nonzero normal 
component, earlier results imposed the Taylor sign condition~\cite{Moran.T.T2015,Moran.T.T2018}, whereas 
Wang--Xin~\cite{Wang.X2024} recently proved local well-posedness without this 
condition, revealing the stabilizing role of the normal magnetic field.

In the present paper, we are interested in the qualitative behavior of such discontinuities. To analyze such high-frequency phenomena, we use the method of nonlinear geometric optics, which is an important and effective approach for studying high-frequency wave propagation in hyperbolic systems. It has been widely used to construct asymptotic solutions in one and several space dimensions \cite{ Hunte.M.R1986, Hunte.K1983, Majda.1986, Majda.A1988, Rauch.2012}. Beyond classical oscillatory wavetrains, nonlinear geometric optics has also
been used to analyze high-frequency phenomena involving boundaries and interfaces. Specifically, it has been applied to the study of surface waves propagating along characteristic boundaries and free interfaces \cite{Ali.H2003, Pierr.C2021, Secch.2015, Secch.Y2022}, as well as to other
resonance and amplification phenomena arising in boundary and interface problems \cite{Kilqu.2022,Coulo.W2017}.

This perspective is particularly relevant for vortex sheets, where the
interaction of oscillatory waves with the interface may produce
resonance and amplification phenomena. Artola and Majda \cite{Artol.M1987,Artol.M1989a,Artol.M1989} studied the nonlinear development of instabilities in supersonic vortex sheets, and identified critical angles of incidence yielding amplification. In a related three-dimensional steady setting, Wang and Yu  established the stability condition of vortex sheets in three-dimensional compressible steady flows and further derived similar amplification phenomena for highly oscillatory vortex sheets~\cite{Wang.Y2013a,Wang.Y2014}.

For the nonisentropic Euler vortex sheets considered here, however, their response to highly oscillatory perturbations remains much less understood. In particular, the possible amplification mechanisms associated with critical angles have not been studied in the same framework.

Motivated by the above works, we study the nonlinear interaction between high-frequency oscillatory waves and a two-dimensional nonisentropic compressible vortex sheet. More precisely, we consider a planar vortex sheet perturbed by a small-amplitude high-frequency incident wave, and analyze the corresponding free boundary problem by means of nonlinear geometric optics. Our purpose is to determine the critical angles of incidence at which amplification occurs, and to derive the effective equations governing the leading profiles of the highly oscillatory vortex sheet. Compared with the isentropic case, the main difficulty comes from the presence of the entropy variable, which changes the characteristic structure of the system, the boundary coupling, and the resulting amplitude equations.
\subsection{Statement of problems}
We recall the full form of the two-dimensional nonisentropic compressible Euler equations:
\begin{equation}
\left\{
\begin{aligned}
    & (\partial_t + \mathbf{u} \cdot \nabla)p + \gamma p \, \nabla \cdot \mathbf{u} = 0, \\
    & \rho (\partial_t + \mathbf{u} \cdot \nabla)\mathbf{u} + \nabla p = 0, \\
    & (\partial_t + \mathbf{u} \cdot \nabla)s = 0.
\end{aligned}
\right.
\label{eq:Euler}
\end{equation}

where $t\ge 0$ and $x=(x_1,x_2)\in\mathbb{R}^2$ denote the time and space variables, respectively. The density is determined by the constitutive law
\begin{equation}
\rho = \rho(p,s) := A p^{\frac{1}{\gamma}} e^{-\frac{s}{\gamma}}, \label{eq:rho}
\end{equation}
where $A>0$ is given and $\gamma>1$ is the adiabatic exponent of the gas.

In this paper, we study perturbations around a piecewise constant vortex sheet, which is set as
\begin{equation}
    \label{constant-state}
    \begin{aligned}
&\Ubar^{R}=(\pbar^{R}, \;  0, \; \ubar_2^{R},\; \sbar^{R}) \quad &&\text{for~} x_1>0,\\
&\Ubar^{L}=(\pbar^{L}, \;  0, \; \ubar_2^{L}, \; \sbar^{L})\quad &&\text{for~} x_1<0,\\
&\text{where~} \ubar_2^R = \bar{v}>0,
\quad
\ubar_2^L = -\bar{v}<0,
\quad
\pbar^{R}=\pbar^{L}=\pbar.
    \end{aligned}
\end{equation}
The superscripts R and L denote functions defined on the right or left of the compressible vortex sheet, respectively. Let us denote by $\bar{c}^{R,L} = c(\bar{p}, \bar{s}^{R,L})$ the sound speeds corresponding to the constant states $\Ubar ^{R,L}$, where for the polytropic gas \eqref{eq:rho},
\begin{equation}
c(p,s) := \sqrt{\frac{\gamma e^{s/\gamma}}{A} p^{1-\frac{1}{\gamma}} }.
\label{eq:c}
\end{equation}

For the two-dimensional nonisentropic Euler equations, the planar vortex sheet is  \emph{linearly stable} if the background states satisfy the stability condition

\begin{equation}
2\bar v \ge \left( (\bar c^{R})^{2/3} + (\bar c^{L})^{2/3} \right)^{3/2}.
\label{eq:linear-stability}
\end{equation}

According to \cite{Morando.2019}, the nonlinear stability condition is
\begin{equation}
2\bar v > \left( (\bar c^{R})^{2/3} + (\bar c^{L})^{2/3} \right)^{3/2},
\qquad
2\bar v \neq \sqrt{2}\,(\bar c^{R}+\bar c^{L}).
\label{eq:nonlinear-stability}
\end{equation}

In what follows, we assume that the background
states satisfy the nonlinear stability condition
\eqref{eq:nonlinear-stability}.
\subsection{Main Result}

We summarize the main conclusions of the formal asymptotic analysis.
\begin{theorem}[Informal Statement]
For \eqref{eq:Euler}, our asymptotic analysis shows that the interaction of small-amplitude, high-frequency oscillatory waves with a nonisentropic vortex sheet leads to resonant amplification at three critical angles of incidence. These critical angles exhibit a twofold asymmetry: the resonant angles associated with waves incident from the same side are not symmetric, and the resonance patterns for waves incident from the two sides of the vortex sheet are different. This nonlinear process is characterized by the emergence of Burgers transport equations  governing the leading-order wave profiles and a Hamilton--Jacobi-type equation describing the evolution of the discontinuity front.
On one resonant branch, the interface equation reduces to a forced linear transport equation. 
\end{theorem}
Precise statements can be found in Theorems~\ref{theorem:resonant-angles} and~\ref{Thm4.2} below.
\begin{remark}
This result extends the classical analysis of resonant amplification, pioneered by Artola and Majda \cite{Artol.M1987} for isentropic vortex sheets, to the more general nonisentropic regime. The primary distinction lies in the presence of the entropy variable. Unlike the isentropic setting, where critical incident angles are typically symmetric, the jump in entropy across the vortex sheet induces an inherent asymmetry in the resonance conditions. This asymmetry appears both in the resonance angles for a fixed incident side and in the distinct resonance patterns generated by waves incident from the two sides. This asymmetry is also reflected in the leading profile equations.
\end{remark}
\subsection{Difficulties and Main Ideas}

The first difficulty comes from the linear degeneracy caused by the nonisentropic Euler system. In contrast to the isentropic case, the normal coefficient matrix in the linearized problem has a zero eigenvalue of multiplicity two~\cite{Morando.2019}. Through analysis, it is found that the zero eigenvalue does not alter the form of the leading wave representation.

The second difficulty concerns the nonlinear geometric optics analysis. The nonisentropic structure changes the acoustic speeds and the boundary coupling on the two sides of the vortex sheet. Consequently, the critical angles differ from those in the isentropic case, and the resonance conditions become intrinsically asymmetric in two respects. For a fixed incident side, the incident, reflected, and refracted angles no longer satisfy the symmetric relations present in the isentropic case. Moreover, since the entropy states differ across the vortex sheet, waves incident from different regions satisfy different resonance conditions and therefore yield different critical angles. This asymmetry is also reflected in the leading profile equations. To handle these effects, we carry out a careful algebraic analysis of the resonance conditions, identify the asymmetric resonant angles, and derive the effective leading-order profile equations at resonance.

\subsection{Outline}
The remainder of this paper is organized as follows. In Sec.~II, we study the linear acoustics for vortex sheets. Sec.~III is devoted to the analysis of the response of a compressible nonisentropic vortex sheet when it is perturbed by high-frequency oscillatory waves, and to the identification of the critical angles of incidence at which the reflected and refracted waves are amplified. In Sec.~IV, we carry out a formal asymptotic expansion for the free-boundary problem of the nonisentropic quasilinear hyperbolic system when a high-frequency oscillatory wave strikes the vortex sheet exactly at a critical angle, and derive the effective equations governing the leading profiles of the highly oscillatory vortex sheets.
\section{Linear acoustics for vortex sheets}
\subsection{The perturbed compressible vortex sheet}

In terms of the conservative variables, the  nonisentropic compressible Euler equations \eqref{eq:Euler} can be written as
\begin{equation}
\frac{\partial {F_0}}{\partial t}+ \frac{\partial{F_1}}{\partial x_1} + \frac{\partial{F_2}}{\partial x_2}  = 0.
\label{eq:Euler conservative}
\end{equation}
In the following study, we use the Taylor expansions of the vector valued function $F_i(U)$, for $i = 0, 1, 2$ as
\begin{equation}
    F_i(U^{R,L} + \varepsilon W) = F_i(U^{R,L}) + \varepsilon A^{R,L}_i W + \frac{\varepsilon^2}{2} B^{R,L}_i (W, W) + \mathcal{O}(\varepsilon^3),
\end{equation}
the Jacobian matrices \(A_i\) $(i=0, 1, 2)$ are defined by 
\begin{equation}
  A_i = \frac{\partial F_i}{\partial w}.
\end{equation}
The expressions of the matrices $A_i$ and $B_i$ are in the appendix; see \eqref{def-A} and \eqref{def-B} below.

We consider the small perturbation of the compressible vortex sheet in the following form
\begin{equation}
    U_\varepsilon(x_1, x_2, t)=\begin{cases}
        U^{R}_\varepsilon=\Ubar^R+ \varepsilon U^{R}_1(x_1, x_2, t, \varepsilon), & x_1>\Phi_{\varepsilon}(x_2,t) \\
        U^{L}_\varepsilon=\Ubar^L+ \varepsilon U^{L}_1(x_1, x_2, t, \varepsilon), & x_1<\Phi_{\varepsilon}(x_2,t)
    \end{cases}
\end{equation}
In order to work in a fixed domain, as in \cite{Artol.M1987}, we introduce new variables to straighten the boundary
\begin{equation}
    \widetilde{x}_1 = x_1- \Phi_\varepsilon(x_2, t), \quad \widetilde{x}_2 = x_2 , \quad \widetilde{t} = t. 
\end{equation}
Under this transformation, the free boundary $\{x_1 =  \Phi_\varepsilon(x_2, t)\}$ is transformed into the fixed one $\{\widetilde{x}_1 = 0\}$, and the domains $\Omega^{R,L} = \{(x_1,x_2, t) \mid \pm (x_1- \Phi_\varepsilon(x_2, t)) > 0\}$ into the half spaces $\widetilde{\Omega}^{R,L} = \{\pm \widetilde{x}_1 > 0\}$.

Setting
\[
\widetilde{U}^{R,L}_\varepsilon(\tilde{x}_1, \tilde{x}_2, \tilde{t}) = U^{R,L}_\varepsilon(x_1, x_2, t), \quad \widetilde{\Phi}_\varepsilon(\tilde{x}_2, \tilde{t}) = \Phi_\varepsilon(x_2, t),
\]
then $(\widetilde{U}^R_\varepsilon, \widetilde{U}^L_\varepsilon, \widetilde{\Phi}_\varepsilon)$ satisfy the following equations (omitting the tilde in the following for simplicity): 
\begin{equation}
\label{eq:F}
\begin{aligned}
    &\partial_t \bigl(F_0(U^{R,L}_\varepsilon)\bigr)
    + \partial_{x_1} \bigl(F_1(U^{R,L}_\varepsilon)\bigr)
    + \partial_{x_2}\bigl(F_2(U^{R,L}_\varepsilon)\bigr) \\
    &\quad
    - \partial_{x_1}\bigl(
        \Phi_{\varepsilon,t} F_0(U^{R,L}_\varepsilon)
        + \Phi_{\varepsilon,x_2}F_2(U^{R,L}_\varepsilon)
    \bigr)
    = 0, \quad
    \begin{aligned}[t]
        &x_1>0 \quad \text{for R},\\
        &x_1<0 \quad \text{for L},
    \end{aligned}
\end{aligned}
\end{equation}

and boundary conditions
\begin{equation}
    \label{eq:boundary}
\begin{cases}
p^{L}_{\varepsilon} - p^{R}_{\varepsilon} = 0, \\[6pt]
\Phi_{\varepsilon,t} + \Phi_{\varepsilon,x_2} 
  u_{2,\varepsilon}^{R,L} - u_{1,\varepsilon}^{R,L} = 0
\end{cases}
\end{equation}

\subsection{Linear geometric optics and causality}
\subsubsection{Very small amplitude sound waves}
We consider asymptotic solutions of the compressible equations consisting of very small amplitude sound waves propagating in all of space through the constant background $\bar{U}^{R,L}$. We seek formal asymptotic solutions with the form
\begin{equation}
    U_\varepsilon^{R,L} = \bar{U}^{R,L} + \varepsilon^2 U_2^{R,L}(x_1,x_2,t, \frac{\phi^{R,L}}{\varepsilon}) + \varepsilon^3 U^{R,L}_3(x_1,x_2,t, \frac{\phi^{R,L}}{\varepsilon}) + O(\varepsilon^4).
    \label{eq:asymptotic expansion} 
\end{equation}
As Artola and Majda \cite{Artol.M1987}, we
assume that \(\phi^{R,L}\) has the simple form
\begin{equation}
    \phi^{R,L} = \omega_0 t + k^{R,L}x_1+  x_2.
    \label{eq:Fast variable}
\end{equation}

By expanding the fluxes, in a standard fashion we have that we can construct a formal solution to the equation within order \(\mathcal{O}(\varepsilon^{2})\) provided that
\begin{equation}
\label{eq:formal solution}
\left(
A_{0}^{R,L}\,\omega_{0}
+ A_{1}^{R,L}k^{R,L}\
+ A_{2}^{R,L}
\right)
\frac{\partial U^{R,L}_{2}}{\partial \theta^{R,L}} = 0.
\end{equation}

Next, from \eqref{eq:Fast variable}:
\begin{equation}
A_{0}^{R,L}\,\omega_{0}
+ A_{1}^{R,L}k^{R,L}\
+ A_{2}^{R,L}
= \mathcal{A}^{R,L}.
\end{equation}
\eqref{eq:formal solution} is equivalent to the following equation
\begin{equation}
\det\!\left( \mathcal{A}^{R,L}\right) = 0.
\end{equation}

\[
\tiny
\mathcal{A}^{R,L}=
\begin{bmatrix}
\dfrac{\omega_0+\bar v^{R,L}}{(\bar c^{R,L})^2}&k^{R,L}\bar\rho^{R,L}&
\bar\rho^{R,L}&-\dfrac{\bar\rho^{R,L}(\omega_0+\bar v^{R,L})}{\gamma}
\\[8pt]
k^{R,L}&\bar\rho^{R,L}(\omega_0+\bar v^{R,L})&0&0
\\[8pt]
1+\dfrac{\bar v^{R,L}(\omega_0+\bar v^{R,L})}{(\bar c^{R,L})^2}
&k^{R,L}\bar\rho^{R,L}\bar v^{R,L}&\bar\rho^{R,L}(\omega_0+2\bar v^{R,L})&
-\dfrac{\bar\rho^{R,L}\bar v^{R,L}(\omega_0+\bar v^{R,L})}{\gamma}
\\[10pt]
\dfrac{\bar s^{R,L}(\omega_0+\bar v^{R,L})}{(\bar c^{R,L})^2}
&k^{R,L}\bar\rho^{R,L}\bar s^{R,L}&\bar\rho^{R,L}\bar s^{R,L}&
\bar\rho^{R,L}(\omega_0+\bar v^{R,L})\!\left(1-\dfrac{\bar s^{R,L}}{\gamma}\right)
\end{bmatrix},
\]
where $\bar v^R=\bar v>0$ and $\bar v^L=-\bar v<0$. Calculating the determinant gives the following result:
\begin{equation}
\det\!\left(\mathcal A^{R,L}\right)=
\frac{\bigl(\bar\rho^{R,L}\bigr)^3}{\bigl(\bar c^{R,L}\bigr)^2}
\bigl(\omega_0+\bar v^{R,L}\bigr)^2
\left[\bigl(\omega_0+\bar v^{R,L}\bigr)^2-
\bigl(\bar c^{R,L}\bigr)^2-\bigl(\bar c^{R,L}\bigr)^2 (k^{R,L})^2
\right]=0.
\end{equation}
where this equation can only solve for two values of \(k^{R,L}\).
\begin{equation}
\label{k}
k^{R,L}=\pm\frac{\sqrt{\bigl(\omega_0+\bar v^{R,L}\bigr)^2-\bigl(\bar c^{R,L}\bigr)^2}
}{\bar c^{R,L}}.
\end{equation}
Notice that in the isentropic case ($\cbar^R=\cbar^L=\cbar$):
\begin{equation}
\label{k in  isentropic}
k^{R,L}=\pm \sqrt{(\frac{\omega_0}{\cbar}+\frac{\bar v^{R,L}}{\cbar})^2-1}
\end{equation}
The results are consistent with those of Majda. We define
\begin{equation}
\label{ a-M}
\alpha^{R,L}=\frac{\omega_0}{\bar c^{R,L}},
\qquad
M^{R,L}=\frac{\bar v^{R,L}}{\bar c^{R,L}}
\end{equation}
Obviously, our determinant has the following solutions:
\begin{equation}
\label{def-kR}
k^{R} = \pm \left( (\alpha^R + M^R)^{2} - 1 \right)^{1/2}, 
\quad |\alpha^R+ M^R| > 1
\end{equation}
Similarly, for the left hand side, we have the following expression:
\begin{equation}
\label{def-kL}
k^{L} = \pm \left( (\alpha^L + M^L)^{2} - 1 \right)^{1/2}, 
\quad |\alpha^L+M^L| > 1.
\end{equation}
Consequently,  \( U_2 \) in \eqref{eq:formal solution} necessarily has the form
\begin{equation}
U_2 = \mathcal{A}(x_1, x_2, t, \theta) r^{R,L}
\label{U_2 with A}
\end{equation}
where \(  \mathcal{A} \)  satisfies the standard linear transport equation of geometric optics.

\subsubsection{Causality}
Now we investigate the propagation directions of the above waves by the causality conditions as \cite{Artol.M1987}.

Based on the ansatz \eqref{eq:asymptotic expansion} and \eqref{U_2 with A}, \eqref{eq:Euler conservative} is satisfied within $\mathcal{O}(\varepsilon^2)$ if and only if the profile $ \mathcal{A}^{R,L}$ satisfies the following linear transport equation:
\begin{equation}
\mathcal{A}_{x_1}^{R,L} 
- \mathcal{A}_{x_2}^{R,L} \frac{\partial k^{R,L}}{\partial \xi} 
- \mathcal{A}_{t}^{R,L} \frac{\partial k^{R,L}}{\partial \omega} = 0,
\label{linear transport equation}
\end{equation}

where
\begin{equation}
k^{R}(\omega,\xi) = \pm \left[ \left(\frac{\omega}{\cbar^R} + M^R\xi\right)^{2} - \xi^{2} \right]^{1/2}, 
\quad
k^{L}(\omega,\xi) = \pm \left[ \left(\frac{\omega}{\cbar^L} + M^L\xi\right)^{2} - \xi^{2} \right]^{1/2},
\end{equation}

\begin{equation}
\frac{\partial k^{R,L}}{\partial \xi} 
= \left. \frac{\partial k^{R,L}(\omega,\xi)}{\partial \xi} \right|_{(\omega_{0},1)}, 
\quad
\frac{\partial k^{R,L}}{\partial \omega} 
= \left. \frac{\partial k^{R,L}(\omega,\xi)}{\partial \omega} \right|_{(\omega_{0},1)},
\end{equation}

\begin{align}
\frac{\partial k^{R,L}}{\partial \xi}\Big|_{(\omega_0,1)} 
&= \pm \frac{M^{R,L}\!\left(\tfrac{\omega_0}{\cbar^{R,L}} + M^{R,L}\right) - 1}
         {\sqrt{\left(\tfrac{\omega_0}{\cbar^{R,L}} + M^{R,L}\right)^2 - 1}}, \\[6pt]
\frac{\partial k^R}{\partial \omega}\Big|_{(\omega_0,1)} 
&= \pm \frac{\tfrac{1}{\cbar^R}\left(\tfrac{\omega_0}{\cbar^R}+M^R\right)}
         {\sqrt{\left(\tfrac{\omega_0}{\cbar^R}+M^R\right)^2 - 1}}, \\[6pt]
\frac{\partial k^L}{\partial \omega}\Big|_{(\omega_0,1)} 
&= \pm \frac{\tfrac{1}{\cbar^L}\left(\tfrac{\omega_0}{\cbar^L}+M^L\right)}
         {\sqrt{\left(\tfrac{\omega_0}{\cbar^L}+M^L\right)^2 - 1}}.
\end{align}

Consequently, we obtain the following results:
\begin{equation}
\begin{aligned}
\text{A)} \quad & \text{Emanates from the free boundary for } x_1>0 \ (x_1<0) \\
& \text{as time increases if } 
\frac{\partial k^R}{\partial \omega} < 0 \ \left( \frac{\partial k^L}{\partial \omega} > 0 \right), \\[10pt]
\text{B)} \quad & \text{Strikes the free boundary for } x_1>0 \ (x_1<0) \\
& \text{as time increases if } 
\frac{\partial k^R}{\partial \omega} > 0 \ \left( \frac{\partial k^L}{\partial \omega} < 0 \right).
\end{aligned}
\end{equation}
In what follows, waves moving away from or toward the contact discontinuity
front \(\{x_1=0\}\) as time increases will be labeled by the subscripts \(+\)
and \(-\), respectively.

\subsubsection{Derivation of the phase function under different conditions}
\label{sect-phase}

To ensure that \eqref{def-kR} and \eqref{def-kL} are well-defined, we have that 
$\alpha^{R,L}$ could be different, i.e., 
\begin{equation}\label{alpha-interval}
\alpha^{R,L} \in
\begin{cases}
(-\infty,-M^R-1)\cup(1-M^R,\infty), & \text{for } \alpha^R,\\
(-\infty,-M^L-1)\cup(1-M^L,\infty), & \text{for } \alpha^L.
\end{cases}
\end{equation}
We will now analyze the variation of the phase in these four cases.

To satisfy the expression in \eqref{linear transport equation} above, we perform a case-by-case analysis for the different intervals, which yields the following results:
\begin{subequations}\label{eq:k-branches}
\begin{align}
&\text{If }\alpha^R\in(-\infty,-M^R-1),\quad
k_+^R = \sqrt{(\alpha^R+M^R)^2-1},\quad
k_-^R = -\sqrt{(\alpha^R+M^R)^2-1},
\label{eq:k-branches-a}
\\
&\text{If }\alpha^R\in(1-M^R,\infty),\quad
k_+^R = -\sqrt{(\alpha^R+M^R)^2-1},\quad
k_-^R = \sqrt{(\alpha^R+M^R)^2-1},
\label{eq:k-branches-b}
\\
&\text{If }\alpha^L\in(-\infty,-M^L-1),\quad
k_+^L = -\sqrt{(\alpha^L+M^L)^2-1},\quad
k_-^L = \sqrt{(\alpha^L+M^L)^2-1},
\label{eq:k-branches-c}
\\
&\text{If }\alpha^L\in(1-M^L,\infty),\quad
k_+^L = \sqrt{(\alpha^L+M^L)^2-1},\quad
k_-^L = -\sqrt{(\alpha^L+M^L)^2-1}.
\label{eq:k-branches-d}
\end{align}
\end{subequations}

The phase functions are expressed as follows (the form varies for different intervals):
\begin{equation}
\begin{gathered}
\phi^{R}_{\pm} = \omega_{0} t + k_{\pm}^{R}x_1 + x_2, \\
\phi^{L}_{\pm} = \omega_{0} t +  k_{\pm}^{L}x_1 +x_2.   
\end{gathered}    
\end{equation}

\subsubsection{Derivation of Eigenvectors under Different Conditions}
Next, using the values of $k$ obtained above, we derive the right and left eigenvectors in two dimensions, respectively:
\begin{equation}
    \label{def-ev-right}
    \left(
A_{0}^{R,L}\,\omega_{0}
+ A_{1}^{R,L}k^{R,L}
+ A_{2}^{R,L}
\right) r^{R,L} = 0,
\end{equation}
\begin{equation}
    \label{def-ev-left}
    l^{R,L} \left(
A_{0}^{R,L}\,\omega_{0}
+ A_{1}^{R,L}k^{R,L}
+ A_{2}^{R,L}
\right) = 0.
\end{equation}

 We obtain a general expression for the eigenvectors, which differs only in the parameter k for different intervals:
\[
r_{\pm}^{R,L} = \left[ -\bar\rho^{R,L}\bigl(\omega_0+\bar v^{R,L}\bigr), \ k_{\pm}^{R,L}, \ 1, \ 0 \right]^\top,
\]
\[
l_{\pm}^{R,L} = \left[ \omega_0+2\bar v^{R,L} -\frac{\bigl(\omega_0+\bar v^{R,L}\bigr)\bar s^{R,L}}{\gamma}, \ -k_{\pm}^{R,L}, \ -1, \ \frac{\omega_0+\bar v^{R,L}}{\gamma} \right]^\top.
\]

\section{The nonlinear response of a compressible vortex sheet to very small amplitude planar sound waves in 2D  }

\subsection{The existence of critical resonant angles of incidence }
We consider a very small amplitude perturbed incident nonlinear sound wave packet
\begin{equation}
\varepsilon^2\mathcal{A}^R_-(x_1,x_2,t,\frac{\phi^R_-}{\varepsilon})r^R_-
\label{incident nonlinear sound wave packet}
\end{equation}

and compute both the response of the compressible vortex sheet as well as the reflected and transmitted sound waves by constructing formal asymptotic solutions of the free-boundary problem. With the above perturbed incident wave \eqref{incident nonlinear sound wave packet}, we do this respecting causality through the ansatz, which is given by:
\begin{equation}
\begin{aligned}
\Phi_{\varepsilon} &= \varepsilon^{3} \Phi_{2}\!\left(x_2,t,\tfrac{\phi_{0}}{\varepsilon}\right) 
+ \mathcal{O}(\varepsilon^{4}), 
\qquad \phi_{0} = \omega_{0} t + x_2, \\[6pt]
U^{R}_{\varepsilon} &=  \bar{U}^R+\varepsilon^{2} \Bigg(
\mathcal{A}^{R}_-\!\left(x_1,x_2,t,\tfrac{\phi_-^{R}}{\varepsilon}\right) r_-^{R} 
+ \mathcal{A}^{R}_{+}\!\left(x_1,x_2,t,\tfrac{\phi^{R}_{+}}{\varepsilon}\right) r^{R}_{+}
\Bigg) + \mathcal{O}(\varepsilon^{3}), 
\qquad x_1 > 0, \\[6pt]
U^{L}_{\varepsilon} &= \bar{U}^L+\varepsilon^{2} \,
\mathcal{A}^{L}_{+}\!\left(x_1,x_2,t,\tfrac{\phi^{L}_{+}}{\varepsilon}\right) r^{L}_{+}
+ \mathcal{O}(\varepsilon^{3}), 
\qquad x_1 < 0 .
\end{aligned}
\label{ansatz}
\end{equation}

We have the following boundary conditions from \eqref{eq:boundary}. The terms of order $\varepsilon^2$ vanish on the boundary provided that:
\begin{equation}
 \mathcal{M}
\begin{pmatrix}
\dfrac{\partial \Phi_{2}}{\partial \theta}(x_2,t,\theta) \\[6pt]
\mathcal{A}^{R}_{+}(0,x_2,t,\theta) \\[6pt]
\mathcal{A}^{L}_{+}(0,x_2,t,\theta)
\end{pmatrix}
=
\mathcal{A}_-^{R}(0,x_2,t,\theta)\,\mathbf{b},
\qquad t>0
\label{eq:linear boundary conditon}
\end{equation}
where the expressions are given by:
\[
\tiny
\mathcal{M} =
\begin{pmatrix}
0 & \rhobar^{R}\bar{c}^R\bigl(\alpha^R+M^{R}\bigr) &-\rhobar^{L}\bar{c}^L\bigl(\alpha^L+M^L\bigr) \\[2ex]
\cbar^R(\alpha^{R}+M^R) &  -k_+^R & 0 \\[2ex]
\cbar^L(\alpha^{L}+M^L) & 0 & -k_+^L \\[2ex]
\end{pmatrix},
\quad
\mathbf{b} =
\begin{pmatrix}
-\rhobar^{R}\bar{c}^R\bigl(\alpha^R+M^{R}\bigr)\\[2ex]
k_-^R\\[2ex]
0
\end{pmatrix}.
\]

Of course, the amplitude $ \mathcal{A}^R_-(x_1,x_2,t,\theta)$ is completely determined from its initial data at t = 0, while causality requires that 
 $ \mathcal{A}^{R,L}_-(x_1,x_2,t,\theta)$ vanish for $t < 0$ .

The asymptotic expansion breaks down when resonance occurs, which corresponds to the condition $\det(\mathcal{M}) = 0$. The determinant is calculated as:
\begin{equation}
\boxed{\det(\mathcal M)
=
\gamma \bar{p}\left[k_+^{L}\bigl(\alpha^R+M^R\bigr)^2
-k_+^{R}\bigl(\alpha^L+M^L\bigr)^2\right]=0}.
\label{Resonant angle}
\end{equation}
Next, we will find the solutions for $\alpha$ that satisfy this resonance condition under different stability regimes. We discuss the general mathematical solutions in equation\eqref{Resonant angle}.

\subsection{Analysis under Stability Conditions}
First, we discuss the interval of existence of the solution. From \eqref{Resonant angle}, $k_+^L$ and $k_+^R$ have the same sign. Moreover, $\alpha^R$ and $\alpha^L$ should have the same sign, their ratio is a positive constant. The domains of existence of $\alpha^R$ and $\alpha^L$ are given in \eqref{eq:k-branches}. \eqref{Resonant angle} can be transformed into the following equation:
\begin{equation}
\boxed{ (U^2 - V^2) \left[ U^2 V^2 - (U^2 + V^2)\right] =0}
\label{eq:detM}
\end{equation}
where
\begin{equation*}
U = \alpha^R + M^R, \quad V = \alpha^L + M^L.
\end{equation*}
Its solution set may be larger than that of \eqref{Resonant angle}, and the admissible roots will be selected from this enlarged set. From nonlinear stability \eqref{eq:nonlinear-stability}:
\begin{equation}
\label{nonlinear-stability1}
2\bar v>\left((\bar c^R)^{2/3}+(\bar c^L)^{2/3}\right)^{3/2}>\bar c^R+\bar c^L
\end{equation}
\begin{equation}
M^R-M^L=\bar v\left(\frac{1}{\bar c^R}+\frac{1}{\bar c^L}\right)
>\frac{\bar c^R+\bar c^L}{2}\left(\frac{1}{\bar c^R}+\frac{1}{\bar c^L}\right)
=\frac{(\bar c^R+\bar c^L)^2}{2\bar c^R\bar c^L}\ge 2.
\end{equation}
From \eqref{Resonant angle},we must have
\begin{equation}
\label{sign condition}
    k_+^Rk_+^L>0
\end{equation}
Since $\alpha^R=\eta\alpha^L$, with $\eta>0$, the two intervals for $\alpha^L$ in \eqref{alpha-interval} can be rewritten in terms of $\alpha^R$. Combined with the two intervals for $\alpha^R$ itself, the admissible range for $\alpha^R$ becomes
\[
\alpha^R\in(-\infty,-M^R-1),\qquad
\alpha^R\in(1-M^R,M^R-\eta\bigr),\qquad
\alpha^R\in\bigl(M^R+\eta,\infty\bigr).
\]
Hence, taking into account the sign condition
\eqref{sign condition}
in \eqref{eq:k-branches} and the nonlinear stability condition \eqref{nonlinear-stability1}, which guarantees the existence of the  middle intervals, we find that the only admissible range is
\begin{equation}
\label{admissible range}
    \boxed{\alpha^R\in(1-M^R,\;M^R-\eta)}
\end{equation}
For $\eta\neq1$, the first factor in \eqref{eq:detM} gives
\begin{equation}
\alpha^R_1 = M^R \frac{1 + \eta}{1 - \eta}, \quad \alpha^R_2 = M^R \frac{1 - \eta}{1 + \eta}
\end{equation}
where
\begin{equation*}
\eta = \frac{\bar{c}^L}{\bar{c}^R}
\end{equation*}

The root $\alpha_1^R$ does not satisfy \eqref{admissible range}, whereas $\alpha_2^R$ is admissible.
When $\eta=1$, directly from $U^2-V^2=4M^R\alpha^R$, the only root of this factor is $\alpha_2^R=0$.
Next, we consider the roots of the second factor in \eqref{eq:detM},
\[
U^2V^2-(U^2+V^2)=0,
\qquad
U=\alpha^R+M^R,\quad
V=\frac{\alpha^R-M^R}{\eta},
\]
for $\alpha^R\in(1-M^R,M^R-\eta)$. Define
\begin{equation}
\label{f}
    f(\alpha^R):=U^2V^2-(U^2+V^2).
\end{equation}

By the sign-change argument in Appendix~\ref{appendix:roots-second-factor}, the quartic $f$ has four distinct real roots under \eqref{nonlinear-stability1}, which we denote by
\begin{equation}
    \alpha^R \in \{\alpha_3^R,\alpha_4^R,\alpha_5^R,\alpha_6^R\}.
\end{equation}
The roots are summarized in the following result.

\begin{lemma}
\label{lemma:roots-order-main}
Assume \eqref{nonlinear-stability1}. Then the equation
\[
U^2V^2-(U^2+V^2)=0,
\qquad
U=\alpha^R+M^R,\quad
V=\frac{\alpha^R-M^R}{\eta},
\]
has four real roots
\[
\alpha_3^R<\alpha_4^R<\alpha_5^R<\alpha_6^R,
\]
and they satisfy
\[
\alpha_3^R<1-M^R<\alpha_4^R<\alpha_5^R<M^R-\eta<\alpha_6^R.
\]

Moreover, set
\[
p:=-(2(M^R)^2+\eta^2+1),\qquad
q:=2M^R(1-\eta^2),\qquad
r:=(M^R)^4-(M^R)^2(\eta^2+1).
\]
Let \(y_*\) be the largest real root of
\[
\bigl((p+y_*)^2-r\bigr)(p+2y_*)=\frac{q^2}{4},
\qquad p+2y_*>0.
\]
Then the two roots in \((1-M^R,M^R-\eta)\) are given by
\[
\alpha_4^R
=
\frac12\left(
-\sqrt{p+2y_*}
+
\sqrt{-\left(3p+2y_*-\frac{2q}{\sqrt{p+2y_*}}\right)}
\right),
\]
and
\[
\alpha_5^R
=
\frac12\left(
\sqrt{p+2y_*}
-
\sqrt{-\left(3p+2y_*+\frac{2q}{\sqrt{p+2y_*}}\right)}
\right).
\]
\end{lemma}

The proof relies only on the explicit computation and analysis of a quartic
equation. We give the details in Appendix~\ref{appendix:roots-second-factor}.
Hence exactly two roots, $\alpha_4^R$ and $\alpha_5^R$, lie in the admissible interval
\[
(1-M^R,\; M^R-\eta).
\]

\begin{theorem}
\label{theorem:resonant-angles}
Assume \eqref{nonlinear-stability1} and
\[
2\bar v\neq \sqrt2\,(\bar c^R+\bar c^L).
\]
Then the determinant equation \eqref{Resonant angle} has exactly three
pairwise distinct real roots in the right admissible range
\eqref{admissible range}, namely
\begin{equation}
\label{resonant angles alpha R}
\alpha^R\in\{\alpha_2^R,\alpha_4^R,\alpha_5^R\},
\end{equation}
where
\[
\alpha_2^R=M^R\frac{1-\eta}{1+\eta},
\]
and \(\alpha_4^R,\alpha_5^R\) are the two roots given in
Lemma~\ref{lemma:roots-order-main}. These are precisely the right-side
incident angles at which the asymptotic expansion breaks down, and hence they
are called the right-side resonant angles.
\end{theorem}
The structure of these three critical angles is closely related to the zeros of the Lopatinskii determinant in the linear stability analysis of nonisentropic vortex sheets~\cite{Moran.T2008}. Similarly, for the left-side incident problem, where the incident wave comes from the left side, the determinant equation has exactly three pairwise distinct real roots in the left admissible range,
\begin{equation}
\label{resonant angles alpha L}
\alpha^L\in\{\alpha_2^L,\alpha_4^L,\alpha_5^L\}.
\end{equation}
The resonant angles are related by
\[
\qquad
\alpha_j^R=\eta\alpha_j^L,
\qquad j=2,4,5.
\]
When $\eta=1$, the two sets of resonant angles coincide. In particular,
$\alpha_2^R=\alpha_2^L=0$. These are the left-side resonant angles.

\begin{remark}
\label{remark3}
The condition
\[
2\bar v\neq \sqrt2\,(\bar c^R+\bar c^L)
\]
excludes the degenerate case where the admissible root
\[
\alpha_2^R=M^R\frac{1-\eta}{1+\eta}
\]
of $U^2-V^2=0$ is simultaneously a root of
\[
U^2V^2-(U^2+V^2)=0.
\]
Indeed, using the definition of \(f\) in \eqref{f}, we have
\[
f(\alpha_2^R)
=
\frac{8(M^R)^2}{(1+\eta)^4}\Bigl(2(M^R)^2-(1+\eta)^2\Bigr),
\]
so $f(\alpha_2^R)=0$ if and only if
\[
2M^R=\sqrt2\,(1+\eta),
\]
which is equivalent to
\[
2\bar v=\sqrt2\,(\bar c^R+\bar c^L).
\]
\end{remark}

\begin{remark}
The characterization in \eqref{resonant angles alpha R} and \eqref{resonant angles alpha L} may be viewed as a
non-isentropic and asymmetric extension of Majda's isentropic result. The ratio
$\eta$ indicates that, due to the asymmetry of the basic states on the two sides of the
contact discontinuity, the critical angles for resonance are in general different
on the right and left sides. Thus, the resonant phenomenon in the present setting
is not merely asymmetric among the critical angles on a fixed side, but also asymmetric between the two sides of the vortex sheet. When $\eta=1$, which corresponds to the isentropic symmetric case, the two sets of resonant angles coincide, namely $\alpha^R=\alpha^L$, and the result reduces to Majda's
classical resonant-angle.
\end{remark}

\section{Small amplitude high-frequency oscillatory perturbations with resonance}
We have analyzed the asymptotic behavior of solutions when the planar contact discontinuity  is perturbed by a small amplitude high frequency oscillatory wave, from which we find there are three specific incident angles of perturbations resulting in the failure of the ansatz given in \eqref{ansatz}. The nonlinear response is incorporated through the ansatz:
\begin{equation}
\scalebox{0.85}{
$\begin{aligned}
\Phi_{\varepsilon} &= \varepsilon^{2} \Phi_{1}\!\left(x_2,t,\tfrac{\phi_{0}}{\varepsilon}\right) +\varepsilon^{3} \Phi_{2}\!\left(x_2,t,\tfrac{\phi_{0}}{\varepsilon}\right) 
+ \mathcal{O}(\varepsilon^{4}), 
\qquad \phi_{0} = \omega_{0} t + x_2, \\[6pt]
U^{R}_{\varepsilon} &=  \bar{U}^R+\varepsilon \,
\mathcal{\sigma}^{R}\!\left(x_1,x_2,t,\tfrac{\phi^{R}_{+}}{\varepsilon}\right) r^{R}_{+}+\varepsilon^{2} \Bigg(
\mathcal{A}^{R}_-\!\left(x_1,x_2,t,\tfrac{\phi_-^{R}}{\varepsilon}\right) r_-^{R} 
+ U^{R}_{2}\!\left(x_1,x_2,t,\tfrac{\phi^{R}_{+}}{\varepsilon},\tfrac{\phi^{R}_{-}}{\varepsilon}\right) 
\Bigg) + \mathcal{O}(\varepsilon^{3}), 
\qquad x_1 > 0, \\[6pt]
U^{L}_{\varepsilon} &= \bar{U}^L+\varepsilon \,
\mathcal{\sigma}^{L}\!\left(x_1,x_2,t,\tfrac{\phi^{L}_{+}}{\varepsilon}\right) r^{L}_{+}+\varepsilon^{2} \,
U^{L}_{2}\!\left(x_1,x_2,t,\tfrac{\phi^{L}_{+}}{\varepsilon},\tfrac{\phi^{L}_{-}}{\varepsilon}\right) 
+ \mathcal{O}(\varepsilon^{3}), 
\qquad x_1 < 0 ,
\end{aligned}$
}
\label{eq:nonlinear ansatz}
\end{equation}
where
\begin{equation}
    \phi_0=\left.\phi_\pm^R\right|_{x_1=0}
=\left.\phi_\pm^L\right|_{x_1=0}.
\end{equation}
Let us derive the problems of the leading profiles $\sigma^{R,L}$ and $\Phi_1$.

\subsection{Second-order asymptotic equation}
We derive the problems of the leading profiles $\sigma^{R,L}$ and $U^{R,L}_2$. First, we consider the solution in {$x_1 > 0$}. We insert the ansatz \eqref{eq:nonlinear ansatz} into the equation of $U^R_\varepsilon$ given in \eqref{eq:F} and collect the terms of order $\varepsilon$.
$\theta^{R,L}_\pm = \frac{\phi^{R,L}_\pm}{\varepsilon}$,
\begin{equation}
    A^R_+ \frac{\partial U^R_2}{\partial \theta^R_+} + A^R_- \frac{\partial U^R_2}{\partial \theta^R_-} = g^R_1(x_1, x_2, t, \theta^R_+) + g^R_2(x_1, x_2, t, \theta^R_+, \theta^R_-), 
\label{U_2}
\end{equation}
where
\begin{align}
    A^R_\pm &= A^R_0 \omega_0+ A^R_1 k^R_\pm + A^R_2, \notag \\
    g^R_1 &= -\left( A^R_0 \frac{\partial \sigma^R}{\partial t} + A^R_1 \frac{\partial \sigma^R}{\partial x_1} + A^R_2 \frac{\partial \sigma^R}{\partial x_2} \right) r^R_+ - \frac{1}{2} \sum_{j=0}^2 \widetilde{\alpha}^R_j B^R_j (r^R_+, r^R_+) \frac{\partial}{\partial \theta^R_+} (\sigma^R)^2,  \label{U_2 terms}\\
    g^R_2 &= (A^R_0 \omega_0 + A^R_2) r^R_+ k^R_+ \frac{\partial \sigma^R}{\partial \theta^R_+} \frac{\partial \Phi_1}{\partial \theta_0} = -(k^R_+)^2 A^R_1 r^R_+  \frac{\partial \sigma^R}{\partial \theta^R_+} \frac{\partial \Phi_1}{\partial \theta_0}(x_2,t,\frac{\theta^R_+ + \theta^R_-}{2}) , \notag
\end{align}
with $\widetilde{\alpha}^R_0 = \omega_0, \widetilde{\alpha}_1^R = k^R_+,   \widetilde{\alpha}_2 ^R= 1$, and $\theta_0 = (k^R_- \theta^R_+ - k^R_+ \theta^R_-)/(k^R_- - k^R_+)$.

To solve $U^R_2$ from \eqref{U_2}, we apply the method given in  \cite{Artol.M1987} to decompose $U^R_2$ as
\begin{equation}
    U^R_2 = U^R_{2,1} + U^R_{2,2}
\end{equation}
such that $U^R_{2,1}$ and $U^R_{2,2}$ satisfy the equations
\begin{equation}
    A^R_+ \frac{\partial U^R_{2,1}}{\partial \theta^R_+} + A^R_- \frac{\partial U^R_{2,1}}{\partial \theta^R_-} = g^R_1(x_1, x_2, t, \theta^R_+)
    \label{U_{2,1}}
\end{equation}
and
\begin{equation}
    A^R_+ \frac{\partial U^R_{2,2}}{\partial \theta^R_+} + A^R_- \frac{\partial U^R_{2,2}}{\partial \theta^R_-} = g^R_2(x_1, x_2, t, \theta^R_+, \theta^R_-)
    \label{U_{2,2}}
\end{equation}
respectively. Define vectors
\begin{equation}\label{vector}
\begin{aligned}
    r_1^{R,L} &= (0,0,1,0)^\top, 
    \quad 
    r_2^{R,L} = (0,0,0,1)^\top, \\
    l_1^{R,L} &= (-\bar v^{R,L},0,1,0),
    \quad
    l_2^{R,L} = (-\bar s^{R,L},0,0,1).
\end{aligned}
\end{equation}

\begin{equation}
    l^R_j A^R_1 = 0,  \quad A^R_1 r^R_j = 0, \quad j = 1, 2.
\end{equation}
It is easy to verify that
\[
    \{r^R_1,r^R_2,r^R_+,r^R_-\}
    \quad\text{and}\quad
    \{l^R_1,l^R_2,l^R_+,l^R_-\}
\]
are two bases of $\mathbb{R}^4$.
Suppose that $U^R_{2,1}(\theta^R_+, \theta^R_-)$ is an almost periodic function satisfying\eqref{U_{2,1}}. Thanks to definition of the left eigenvectors\eqref{def-ev-left}, multiplying \eqref{U_{2,1}} by $l^R_+$ gives
\begin{equation}
    l^R_+ A^R_- \frac{\partial U^R_{2,1}}{\partial \theta^R_-} = l^R_+ g^R_1(x_1, x_2, t,\theta^R_+).
    \label{multiplying l}
\end{equation}
Since $l^R_+ A^R_- \neq 0$, and $l^R_+ g^R_1(\theta^R_+)$ is independent of $\theta^R_-$, it follows from \eqref{multiplying l} that
\begin{equation} 
    l^R_+ g^R_1(x_1, x_2, t,\theta^R_+) = 0. 
\label{independent of theta}
\end{equation}
By using a direct computation from the expressions of the right and left eigenvectors, one has
\begin{equation}
    \frac{l_+^R A_2^R r_+^R}{l_+^R A_1^R r_+^R}
    =
    -\frac{\partial k_+^R}{\partial \xi},
    \qquad
    \frac{l_+^R A_0^R r_+^R}{l_+^R A_1^R r_+^R}
    =
    -\frac{\partial k_+^R}{\partial \omega},
    \label{eq:transport-coefficients-R}
\end{equation}
where \(l_+^R A_1^R r_+^R\neq 0\).  Combining \eqref{eq:transport-coefficients-R} with the expression of \(g_1^R\) given in \eqref{U_2 terms}, it follows that \(\sigma^R(x_1,x_2,t,\theta^R_+)\) satisfies the following nonlinear Burgers-transport equation in \(\{x_1>0\}\):
\begin{equation} \label{eq:Burgers-sigma}
\begin{aligned}
    & \frac{\partial \sigma^R}{\partial x_1} - \frac{\partial k^R_+}{\partial \xi} \frac{\partial  \sigma^R}{\partial x_2} - \frac{\partial k^R_+}{\partial \omega} \frac{\partial  \sigma^R}{\partial t} + \tilde{b}^R \frac{\partial}{\partial \theta^R_+} ( \sigma^R)^2 = 0, \quad t> 0, \\
    &  \sigma^R(x_1,x_2,t, \theta^R_+) = 0, \quad t \leq 0,
\end{aligned} 
\end{equation}
where
\begin{align*}
    \tilde{b}^R &= \frac{1}{2} (l^R_+ A^R_1 r^R_+)^{-1} l^R_+ \sum_{j=0}^2\widetilde{\alpha}^R_j B^R_j (r^R_+, r^R_+)
\end{align*}

We are going to derive the representation of the solution $U_2^R$ to \eqref{U_2}, providing that $\sigma^R(x_1, x_2, t, \theta_+^R)$ and $\Phi_1(x_2, t, \theta_0)$ are determined already. With the condition in \eqref{independent of theta}, we seek a solution with the form 
\begin{equation}
     \label{eq:U_21}
     U^R_{2,1}=a^R_1(\theta_+^R)r^R_1+a^R_2(\theta_+^R)r^R_2+a^R_3(\theta_+^R,\theta_-^R)r_-^R
\end{equation} 
By a simple computation, we have
\begin{equation*}
    l_-^RA_+^Rr_1^R=l_-^RA_+^Rr_2^R=0
\end{equation*}
so by multiplying $l_-^R$ to \eqref{U_{2,1}} and using the representation \eqref{eq:U_21}, one gets
\begin{equation}
    a_3^R=(l_-^RA_1^Rr_-^R)^{-1}(k_+^R-k_-^R)^{-1}(\int_{\theta_-^R}^{\theta_+^R}l_-^Rg^R_1(x_1, x_2, t, \theta)d\theta)
\end{equation}
Similarly, we have
\begin{equation*}
l_1^R A_\pm^R r_2^R = 0,\qquad
l_1^R A_\pm^R r_-^R = 0,\qquad
l_1^R A_\pm^R r_1^R = \bar\rho^R(\omega_0+\bar v^R),
\end{equation*}
\begin{equation*}
l_2^R A_\pm^R r_2^R = \bar\rho^R(\omega_0+\bar v^R),\qquad
l_2^R A_\pm^R r_-^R = 0,\qquad
l_2^R A_\pm^R r_1^R = 0.
\end{equation*}
Hence,
\begin{equation}
\bar\rho^R(\omega_0+\bar v^R)
\bigl(\partial_{\theta_+^R} a_2^R+\partial_{\theta_-^R} a_2^R\bigr)
=l_2^R g^R_1(x_1,x_2,t,\theta_+^R),
\end{equation}
\begin{equation}
\bar\rho^R(\omega_0+\bar v^R)
\bigl(\partial_{\theta_+^R} a_1^R+\partial_{\theta_-^R} a_1^R\bigr)
=l_1^R g^R_1(x_1,x_2,t,\theta_+^R),
\end{equation}
which implies
\begin{equation}
   a_1^R=(\bar\rho^R(\omega_0+\bar v^R))^{-1}\int_0^{\theta^R_+}l_1^Rg^R_1(x_1, x_2, t, \theta)d\theta,
\end{equation}
\begin{equation}
   a_2^R=(\bar\rho^R(\omega_0+\bar v^R))^{-1}\int_0^{\theta^R_+}l_2^Rg^R_1(x_1, x_2, t, \theta)d\theta,
\end{equation}
being independent of $\theta_-^R$. And we construct a solution to \eqref{U_{2,2}} with the form 
\begin{equation}
\label{The expansion of U_22}
     U^R_{2,2}=a^R_4 (\theta_+^R,\theta_-^R)r_+^R
\end{equation}
so by multiplying $l_+^R$ to \eqref{U_{2,2}} and using the representation \eqref{The expansion of U_22}, one gets
\begin{equation}
    a_4^R=(l_+^RA_1^Rr_+^R)^{-1}(k_-^R-k_+^R)^{-1}(\int_{\theta_+^R}^{\theta_-^R}l_+^Rg^R_2(x_1, x_2, t, \theta)d\theta)
\end{equation}
We next derive the representation of the second-order profile in the region
$\{x_1<0\}$. Inserting the ansatz \eqref{eq:nonlinear ansatz} into the
equation of $U^L_\varepsilon$ and collecting the terms of order $\varepsilon$,
we obtain
\begin{equation}
    A^L_+ \frac{\partial U^L_2}{\partial \theta^L_+}
    +
    A^L_- \frac{\partial U^L_2}{\partial \theta^L_-}
    =
    g^L_1(x_1,x_2,t,\theta^L_+)
    +
    g^L_2(x_1,x_2,t,\theta^L_+,\theta^L_-),
\label{U_2_L}
\end{equation}
where
\begin{align}
    A^L_\pm
    &=
    A^L_0\omega_0+A^L_1 k^L_\pm+A^L_2,
    \notag\\
    g^L_1
    &=
    -\left(
    A^L_0\frac{\partial \sigma^L}{\partial t}
    +
    A^L_1\frac{\partial \sigma^L}{\partial x_1}
    +
    A^L_2\frac{\partial \sigma^L}{\partial x_2}
    \right)r^L_+
    -
    \frac12
    \sum_{j=0}^2
    \widetilde{\alpha}^L_j B^L_j(r^L_+,r^L_+)
    \frac{\partial}{\partial \theta^L_+}(\sigma^L)^2,
    \label{U_2_terms_L}\\
    g^L_2
    &=
    (A^L_0\omega_0+A^L_2)r^L_+ k^L_+
    \frac{\partial \sigma^L}{\partial \theta^L_+}
    \frac{\partial \Phi_1}{\partial \theta_0}=
    -(k^L_+)^2 A^L_1 r^L_+
    \frac{\partial \sigma^L}{\partial \theta^L_+}
    \frac{\partial \Phi_1}{\partial \theta_0}
    \left(
    x_2,t,\frac{\theta^L_+ + \theta^L_-}{2}
    \right).
    \notag
\end{align}
Here
\[
   \widetilde{\alpha}^L_0=\omega_0,
    \widetilde{\alpha}^L_1=k^L_+,
   \widetilde{\alpha}^L_2=1,
    \theta_0
    =
    \frac{k^L_- \theta^L_+ - k^L_+ \theta^L_-}
    {k^L_- - k^L_+}.
\]

Following the same steps, we obtain the
Burgers-transport equation for $\sigma^L$ in $\{x_1<0\}$:
\begin{equation}
\label{eq:Burgers-sigma-L}
\begin{aligned}
    &
    \frac{\partial \sigma^L}{\partial x_1}
    -
    \frac{\partial k^L_+}{\partial \xi}
    \frac{\partial \sigma^L}{\partial x_2}
    -
    \frac{\partial k^L_+}{\partial \omega}
    \frac{\partial \sigma^L}{\partial t}
    +
    \tilde b^L
    \frac{\partial}{\partial \theta^L_+}(\sigma^L)^2
    =
    0,
    \qquad t>0,
    \\
    &
    \sigma^L(x_1,x_2,t,\theta^L_+)=0,
    \qquad t\leq 0,
\end{aligned}
\end{equation}
where
\begin{equation}
    \tilde b^L
    =
    \frac12
    \left(l^L_+ A^L_1 r^L_+\right)^{-1}
    l^L_+
    \sum_{j=0}^2
    \widetilde{\alpha}^L_j B^L_j(r^L_+,r^L_+).
\end{equation}

From \eqref{vector}, a direct computation gives
\begin{equation}
    l^L_j A^L_1=0,\qquad
    A^L_1 r^L_j=0,\qquad j=1,2.
\end{equation}
Similarly, one has
\[
    \{r^L_1,r^L_2,r^L_+,r^L_-\}
    \quad\text{and}\quad
    \{l^L_1,l^L_2,l^L_+,l^L_-\}
\]
are two bases of $\mathbb{R}^4$.
By a similar direct computation as in the right-hand side case, we have
\begin{equation*}
    l^L_- A^L_+ r^L_1=0,\qquad
    l^L_- A^L_+ r^L_2=0.
\end{equation*}
Moreover,
\begin{equation*}
    l^L_1 A^L_\pm r^L_2=0,\qquad
    l^L_1 A^L_\pm r^L_-=0,\qquad
    l^L_1 A^L_\pm r^L_1
    =
    \bar\rho^L(\omega_0+\bar v^L),
\end{equation*}
and
\begin{equation*}
    l^L_2 A^L_\pm r^L_2
    =
    \bar\rho^L(\omega_0+\bar v^L),\qquad
    l^L_2 A^L_\pm r^L_-=0,\qquad
    l^L_2 A^L_\pm r^L_1=0.
\end{equation*}

In a way similar to the derivation from \eqref{eq:U_21} to \eqref{The expansion of U_22} for determining $U^R_2$, we finally obtain the following representation.
\begin{lemma}
\label{lemma:U2L_representation}
The profile $U^L_2$ of \eqref{U_2_L} is given by
\begin{equation}
    U^L_2
    =
    a^L_1 (\theta_+^L)r^L_1
    +
    a^L_2  (\theta_+^L)r^L_2
    +
    a^L_3  (\theta_+^L,\theta_-^L)r^L_-
    +
    a^L_4 (\theta_+^L,\theta_-^L) r^L_+,
\end{equation}
where
\begin{align}
   a^L_1
   &=
   \left(\bar\rho^L(\omega_0+\bar v^L)\right)^{-1}
   \int_0^{\theta^L_+}
   l^L_1 g^L_1(x_1,x_2,t,\theta)\,d\theta,
   \\
   a^L_2
   &=
   \left(\bar\rho^L(\omega_0+\bar v^L)\right)^{-1}
   \int_0^{\theta^L_+}
   l^L_2 g^L_1(x_1,x_2,t,\theta)\,d\theta,
   \\
   a^L_3
   &=
   \left(l^L_- A^L_1 r^L_-\right)^{-1}
   \left(k^L_+-k^L_-\right)^{-1}
   \int_{\theta^L_-}^{\theta^L_+}
   l^L_- g^L_1(x_1,x_2,t,\theta)\,d\theta,
   \\
   a^L_4
   &=
   \left(l^L_+ A^L_1 r^L_+\right)^{-1}
   \left(k^L_--k^L_+\right)^{-1}
   \int_{\theta^L_+}^{\theta^L_-}
   l^L_+ g^L_2(x_1,x_2,t,\theta)\,d\theta.
\end{align}
\end{lemma}

\subsection{The boundary conditions to second order at the perturbed slip-stream surface}
Considering the boundary conditions given in \eqref{eq:boundary}, using the ansatz \eqref{eq:nonlinear ansatz} of $U_{R,L}^\varepsilon$ and collecting the terms of order $\varepsilon$, it yields that
\begin{equation}
\mathcal{M}(\alpha_\mu)
\begin{pmatrix}
\dfrac{\partial \Phi_{1}}{\partial \theta_0}(x_2,t,\theta_0) \\[6pt]
\sigma^R (0,x_2,t,\theta_0)\\[6pt]
\sigma^L (0,x_2,t,\theta_0)
\end{pmatrix}
=0    
\label{eq:boundary 1st}
\end{equation}
where the matrix $\mathcal{M}(\alpha^{R,L}_\mu)$ is defined in \eqref{eq:linear boundary conditon}. 
It has been shown that $\det \mathcal{M}(\alpha^{R,L}_\mu) = 0$ if $\alpha^{R,L}_\mu$ belongs to the set given in \eqref{resonant angles alpha R} and  \eqref{resonant angles alpha L} . It is easy to check from the eigenvectors that
\eqref{eq:boundary 1st} is satisfied if and only if the following two constraints are enforced:
\begin{equation}
\begin{aligned}
   \sigma^R (0,x_2,t,\theta_0) &= \frac{\cbar^R(\alpha^{R}_\mu+M^R)}{k_+^R}\frac{\partial \Phi_{1}}{\partial \theta_0}(x_2,t,\theta_0), \\
   \sigma^L (0,x_2,t,\theta_0) &= \frac{\cbar^L(\alpha^{L}_\mu+M^L)}{k_+^L}\frac{\partial \Phi_{1}}{\partial \theta_0}(x_2,t,\theta_0).
\end{aligned}
\label{eq:Burgers-boundary}
\end{equation}
Restricting \eqref{eq:Burgers-sigma} on $x_1=0$, with \eqref{eq:Burgers-boundary}, we will obtain the Burgers-transport equation for $\Phi_1$.  
Next, we include nonlinear boundary effects at second order and derive a Hamilton-Jacobi equation for $\Phi_1$. We define $U_2^{R,L}=U^{R,L}_{2,1}+U^{R,L}_{2,2}+U^{R,L}_{2,3}$, where 
\begin{equation}
    U^{R,L}_{2,3}=\sigma^{R,L}_2(\theta_+)r^{R,L}_+
\end{equation}
It can be simply verified that adding $U^{R,L}_{2,3}$ does not change the formula \eqref{U_2}. With the ansatz in \eqref{eq:nonlinear ansatz}, we expand the free boundary conditions at the transformed free boundary to second order in $\varepsilon$; the terms of order $\varepsilon^2$ vanish if and only if
\begin{equation}
 \mathcal{M}(\alpha^{R,L}_\mu)
\begin{pmatrix}
\dfrac{\partial \Phi_{2}}{\partial \theta}(x_2,t,\theta_0) \\[6pt]
\sigma_2^R(\theta_+^R)\\[6pt]
\sigma_2^L(\theta_+^L)
\end{pmatrix}
=\begin{pmatrix}
0 \\[6pt]
-(\dfrac{\partial \Phi_{1}}{\partial t}+\bar{c}^RM^R\dfrac{\partial \Phi_{1}}{\partial x_2}+\dfrac{\partial \Phi_{1}}{\partial \theta_0} \sigma^R)\\[6pt]
-(\dfrac{\partial \Phi_{1}}{\partial t}+\bar{c}^LM^L\dfrac{\partial \Phi_{1}}{\partial x_2}+\dfrac{\partial \Phi_{1}}{\partial \theta_0} \sigma^L)\\
\end{pmatrix}+
\mathcal{A}_-^{R}(0,x_2,t,\theta)\mathbf{b}
\qquad t>0
\label{eq:second order boundary}
\end{equation}
\[
\mathbf{b} =\begin{pmatrix}
-\rhobar^{R}\bar{c}^R\bigl(\alpha_\mu^R+M^{R}\bigr)\\[6pt]
k_-^R\\[6pt]
0
\end{pmatrix},\quad \sigma^{R,L}=  \sigma^{R,L} (0,x_2,t,\theta_0)
\]
 In the computation determining the equation in \eqref{eq:second order boundary}, we have used the information that both $U^{R,L}_{2,1}$ and $U^{R,L}_{2,2}$ make no contributions to the boundary at second order. The first reason is that at the boundary $x_1=0$, we have $\theta^{R,L}_+ = \theta^{R,L}_- = \theta_0$, which implies that the coefficients of $U^{R,L}_2$ (except for that of $r_1$ and $r_2$) vanish. The second reason concerns the coefficient of $r_1$ and $r_2$; upon substitution into the boundary conditions, it is evident that this term appears only at the order of $O(\varepsilon^3)$. Thus, it is absent from  \eqref{eq:second order boundary}. By a simple computation, the left eigenvector $l_\mu$ satisfying $l_\mu \mathcal{M}(\alpha_\mu)=0$ is given by
\begin{equation}
    l_\mu=\left(1,\,\frac{\rhobar^R\bar c^R(\alpha_\mu^R+M^R)}{k_+^R},\,-\frac{\rhobar^L\bar c^L(\alpha_\mu^L+M^L)}{k_+^L}\right).
\end{equation}
Multiplying  \eqref{eq:second order boundary} by $l_\mu$, we obtain the following equation for $\Phi_1$
\begin{equation}
    a_\mu \frac{\partial \Phi_{1}}{\partial t} + b_\mu \frac{\partial \Phi_{1}}{\partial x_2} + c_\mu \left(\frac{\partial \Phi_{1}}{\partial \theta_0}\right)^2 = d_\mu \mathcal{A}_-^{R}(0,x_2,t,\theta_0),
    \label{eq:Hamilton-Jacobi equation}
\end{equation}
where the coefficients are defined by the components of the left eigenvector $l_\mu=(1, l_2, l_3)$ as follows:
\begin{equation}\label{H-J-components}
\begin{aligned}
a_\mu &= l_2+l_3,\\
b_\mu &= \bar c^R l_2 M^R+\bar c^L l_3 M^L,\\
c_\mu &= l_2 \frac{\bar c^R(\alpha_\mu^R+M^R)}{k_+^R}
      + l_3 \frac{\bar c^L(\alpha_\mu^L+M^L)}{k_+^L},\\
d_\mu &= -\,l_2 (k_+^R-k_-^R).
\end{aligned}
\end{equation}

Using $\bar\rho^{R,L}(\bar c^{R,L})^2=\gamma\bar p$ in \eqref{H-J-components}, we obtain
\[
c_\mu=\gamma\bar p
\left(\frac{U^2}{U^2-1}-\frac{V^2}{V^2-1}\right),
\quad
U=\frac{\omega_0+\bar v}{\bar c^R},\quad
V=\frac{\omega_0-\bar v}{\bar c^L}.
\]
At $\alpha_2^R=M^R(1-\eta)/(1+\eta)$, we have $U=-V$ and hence $c_2=0$.
On this branch, \eqref{eq:Hamilton-Jacobi equation} reduces to a forced linear transport equation for $\Phi_1$, while the bulk Burgers--transport equations remain nonlinear.

At each resonant angle $\alpha_\mu$, we have $a_\mu\neq 0$.
Thus, $t$ is the evolution variable for the formal leading-order equation \eqref{eq:Hamilton-Jacobi equation}.

\begin{theorem}\label{Thm4.2}

When a high-frequency incident wave strikes the front at a resonant parameter
$\alpha^R=\alpha_\mu^R$ from \eqref{resonant angles alpha R}, the formal
asymptotic expansion \eqref{eq:nonlinear ansatz} of \eqref{eq:F}--\eqref{eq:boundary}
yields the following leading-profile equations, as long as the profiles remain
sufficiently smooth. The initial profile data are assumed
to satisfy the initial--boundary compatibility
conditions associated with \eqref{eq:Burgers-boundary}:

\begin{enumerate}
\item The profile $\Phi_{1}(x_2,t,\theta_0)$ satisfies the Cauchy problem of the Hamilton-Jacobi equation
\[
\left\{
\begin{array}{ll}
\displaystyle
 a_\mu \frac{\partial \Phi_{1}}{\partial t}
 + b_\mu \frac{\partial \Phi_{1}}{\partial x_2}
 + c_\mu \left(\frac{\partial \Phi_{1}}{\partial \theta_0}\right)^2
 = d_\mu \mathcal{A}_-^{R}(0,x_2,t,\theta_0),
& t>0, \\[0.8em]
\displaystyle
 \Phi_1(x_2,0,\theta_0)=\Phi_1^0(x_2,\theta_0).
&
\end{array}
\right.
\]

where $a_{\mu}$, $b_{\mu}$, $c_{\mu}$, and $d_{\mu}$ are given in \eqref{H-J-components}.
In particular, when $\mu=2$, one has $c_2=0$, so the interface equation
reduces to a forced linear transport equation.

\item The profile $\sigma^R(x_1,x_2,t,\theta^R_+)$ solves the initial-boundary value problem of a Burgers-transport equation
\[
 \frac{\partial \sigma^R}{\partial x_1} - \frac{\partial k^R_+}{\partial \xi} \frac{\partial  \sigma^R}{\partial x_2} - \frac{\partial k^R_+}{\partial \omega} \frac{\partial  \sigma^R}{\partial t} + \tilde{b}^R \frac{\partial}{\partial \theta^R_+} ( \sigma^R)^2 = 0, \quad t> 0,x_1>0
\]
with initial-boundary conditions
\[
\sigma^R|_{t\le 0} = 0, \quad \sigma^R|_{x_1=0} =\frac{\cbar^R(\alpha^{R}_\mu+M^R)}{k_+^R}\frac{\partial \Phi_{1}}{\partial \theta_0}(x_2,t,\theta_0)
\]
and
\[
\tilde{b}^R = \frac{1}{2} (l^R_+ A^R_1 r^R_+)^{-1} l^R_+ \sum_{j=0}^2 \widetilde{\alpha}^R_j B^R_j (r^R_+, r^R_+)
\]

\item Similarly, the profile $\sigma^L(x_1,x_2,t,\theta^L_+)$ solves the initial-boundary value problem of a Burgers-transport equation
\[
\frac{\partial \sigma^L}{\partial x_1}
    -
    \frac{\partial k^L_+}{\partial \xi}
    \frac{\partial \sigma^L}{\partial x_2}
    -
    \frac{\partial k^L_+}{\partial \omega}
    \frac{\partial \sigma^L}{\partial t}
    +
    \tilde b^L
    \frac{\partial}{\partial \theta^L_+}(\sigma^L)^2
    =
    0,
    \qquad t>0,\quad x_1<0
\]
with initial-boundary conditions
\[
\sigma^L|_{t\le 0} = 0, \quad \sigma^L|_{x_1=0} =\frac{\cbar^L(\alpha^{L}_\mu+M^L)}{k_+^L}\frac{\partial \Phi_{1}}{\partial \theta_0}(x_2,t,\theta_0),
\]
and
\[
 \tilde b^L
    =
    \frac12
    \left(l^L_+ A^L_1 r^L_+\right)^{-1}
    l^L_+
    \sum_{j=0}^2
   \widetilde{\alpha}^L_j B^L_j(r^L_+,r^L_+).
\]
\end{enumerate}
\end{theorem}
\begin{remark}
The preceding formal asymptotic analysis predicts resonant amplification of
high-frequency oscillatory perturbations. When an incoming wave of amplitude
\(O(\varepsilon^2)\) impinges on the contact discontinuity at a resonant angle,
the generated reflected and refracted waves have amplitudes of order
\(O(\varepsilon)\). Thus the generated outgoing waves are amplified by one order in
\(\varepsilon^{-1}\).
This resonant amplification may be viewed as a physical indication of instability.

This phenomenon is different from the isentropic symmetric case in two respects.
First, the resonant angles on the two sides of the discontinuity are generally
not the same, because the thermodynamic states and sound speeds are different
across the front. Second, the leading profiles on the two sides are governed by
intrinsically distinct Burgers--transport equations, in contrast with the isentropic setting where the two-sided Burgers dynamics are mirror images of a single equation.
The same intrinsic asymmetry also propagates into the Hamilton–Jacobi equation. Therefore, in the present non-isentropic setting, the resonance is not only an
amplification mechanism, but also reflects the intrinsic asymmetry
of the underlying basic states.
\end{remark}
\section*{ACKNOWLEDGMENTS}

This work is partially supported by the National Natural Science Foundation of China (No.12471210). The authors would like to thank Prof. Yuan Yuan for his insightful comments and invaluable guidance throughout this study. 
\section*{AUTHOR DECLARATIONS}

\subsection*{Conflict of Interest}
The authors have no conflicts to disclose.

\subsection*{Author Contributions}
Zihao Li: Conceptualization (equal); Formal analysis (equal); Investigation; Writing -- original draft; Writing -- review \& editing (equal).

Yang Xu: Conceptualization (equal); Supervision; Formal analysis (equal); Writing -- review \& editing (equal).
\section*{DATA AVAILABILITY}
Data sharing is not applicable to this article as no new data were created or analyzed in this study.
\appendix
\section{Matrices}
The vectors \(F_0, F_1, F_2\) are given by:
\[
  w=\begin{bmatrix}
     p\\
     u_1 \\
     u_2 \\
    s
  \end{bmatrix},
  \qquad
  F_0 =\begin{bmatrix}
     \rho\\
    \rho u_1 \\
    \rho u_2 \\
    \rho s
  \end{bmatrix},
  \qquad
  F_1 =
  \begin{bmatrix}
    \rho u_1\\
    \rho u_1^2 + p \\
    \rho u_1 u_2 \\
    \rho s u_1
    
  \end{bmatrix},
  \qquad
  F_2=
   \begin{bmatrix}
    \rho u_2\\
    \rho u_1 u_2  \\
    \rho u_2^2  +p\\
    \rho s u_2
   \end{bmatrix}
    \qquad 
\]
The Jacobian matrices $A_i(W)$ and $B_i(W,W)$ are as follows. For simplicity, we omit the superscripts $R,L$ below.
{\tiny
\begin{align}
A_0&=
\begin{bmatrix}
\dfrac{1}{(\bar{c})^2} & 0 & 0 & -\dfrac{\bar\rho}{\gamma}\\[8pt]
0 & \bar\rho & 0 & 0\\[4pt]
\dfrac{\bar v}{(\bar{c})^2} & 0 & \bar\rho & -\dfrac{\bar\rho\bar v}{\gamma}\\[8pt]
\dfrac{\bar s}{(\bar{c})^2} & 0 & 0 & \bar\rho\!\left(1-\dfrac{\bar s}{\gamma}\right)
\end{bmatrix}, \notag\\[10pt]
A_1&=
\begin{bmatrix}
0 & \bar\rho & 0 & 0\\
1 & 0 & 0 & 0\\
0 & \bar\rho\bar v & 0 & 0\\
0 & \bar\rho\bar s & 0 & 0
\end{bmatrix}, \label{def-A}\\[10pt]
A_2&=
\begin{bmatrix}
\dfrac{\bar v}{(\bar{c})^2} & 0 & \bar\rho & -\dfrac{\bar\rho\bar v}{\gamma}\\[8pt]
0 & \bar\rho\bar v & 0 & 0\\[4pt]
1+\dfrac{(\bar v)^2}{(\bar{c})^2} & 0 & 2\bar\rho\bar v & -\dfrac{\bar\rho(\bar v)^2}{\gamma}\\[8pt]
\dfrac{\bar s\bar v}{(\bar{c})^2} & 0 & \bar\rho\bar s &
\bar\rho\bar v\!\left(1-\dfrac{\bar s}{\gamma}\right)
\end{bmatrix}. \notag
\end{align}
}
{\tiny
\begin{align}
B_0(W,W) &=
\begin{bmatrix}
-\dfrac{(\gamma-1)}{\gamma\bar p(\bar c)^2}W_1^2
-\dfrac{2}{\gamma(\bar c)^2}W_1W_4
+\dfrac{\bar p}{\gamma(\bar c)^2}W_4^2
\\[12pt]
\dfrac{2}{(\bar c)^2}W_1W_2
-\dfrac{2\bar p}{(\bar c)^2}W_2W_4
\\[12pt]
-\dfrac{(\gamma-1)\bar v}{\gamma\bar p(\bar c)^2}W_1^2
-\dfrac{2\bar v}{\gamma(\bar c)^2}W_1W_4
+\dfrac{\bar p\,\bar v}{\gamma(\bar c)^2}W_4^2
+\dfrac{2}{(\bar c)^2}W_1W_3
-\dfrac{2\bar p}{(\bar c)^2}W_3W_4
\\[12pt]
-\dfrac{(\gamma-1)\bar s}{\gamma\bar p(\bar c)^2}W_1^2
+\dfrac{2\left(1-\bar s/\gamma\right)}{(\bar c)^2}W_1W_4
+\dfrac{\bar p\left(\bar s/\gamma-2\right)}{(\bar c)^2}W_4^2
\end{bmatrix}, \notag \\[18pt]
B_1(W,W) &=
\begin{bmatrix}
\dfrac{2}{(\bar c)^2}W_1W_2
-\dfrac{2\bar p}{(\bar c)^2}W_2W_4
\\[12pt]
2\bar\rho W_2^2
\\[12pt]
\dfrac{2\bar v}{(\bar c)^2}W_1W_2
+2\bar\rho W_2W_3
-\dfrac{2\bar p\,\bar v}{(\bar c)^2}W_2W_4
\\[12pt]
\dfrac{2\bar s}{(\bar c)^2}W_1W_2
+2\left(\bar\rho-\dfrac{\bar p\,\bar s}{(\bar c)^2}\right)W_2W_4
\end{bmatrix}, \notag \\[18pt]
B_2(W,W) &=
\begin{bmatrix}
-\dfrac{(\gamma-1)\bar v}{\gamma\bar p(\bar c)^2}W_1^2
-\dfrac{2\bar v}{\gamma(\bar c)^2}W_1W_4
+\dfrac{\bar p\,\bar v}{\gamma(\bar c)^2}W_4^2
+\dfrac{2}{(\bar c)^2}W_1W_3
-\dfrac{2\bar p}{(\bar c)^2}W_3W_4
\\[12pt]
\dfrac{2\bar v}{(\bar c)^2}W_1W_2
+2\bar\rho W_2W_3
-\dfrac{2\bar p\,\bar v}{(\bar c)^2}W_2W_4
\\[12pt]
-\dfrac{(\gamma-1)(\bar v)^2}{\gamma\bar p(\bar c)^2}W_1^2
-\dfrac{2(\bar v)^2}{\gamma(\bar c)^2}W_1W_4
+\dfrac{\bar p(\bar v)^2}{\gamma(\bar c)^2}W_4^2
+\dfrac{4\bar v}{(\bar c)^2}W_1W_3
+2\bar\rho W_3^2
-\dfrac{4\bar p\,\bar v}{(\bar c)^2}W_3W_4
\\[12pt]
\begin{aligned} 
&-\dfrac{(\gamma-1)\bar v\bar s}{\gamma\bar p(\bar c)^2}W_1^2
+\dfrac{2\bar v\left(1-\bar s/\gamma\right)}{(\bar c)^2}W_1W_4
+\dfrac{\bar p\,\bar v\left(\bar s/\gamma-2\right)}{(\bar c)^2}W_4^2 +\dfrac{2\bar s}{(\bar c)^2}W_1W_3
+2\left(\bar\rho-\dfrac{\bar p\,\bar s}{(\bar c)^2}\right)W_3W_4
\end{aligned}
\end{bmatrix}. \label{def-B}
\end{align}
}

\section{Proof of Lemma~\ref{lemma:roots-order-main}}
\label{appendix:roots-second-factor}
\begin{proof}
Let
\[
f(\alpha^R):=U^2V^2-(U^2+V^2),
\qquad
U=\alpha^R+M^R,\quad
V=\frac{\alpha^R-M^R}{\eta}.
\]
Then $f$ is a quartic polynomial in $\alpha^R$ with positive leading coefficient, hence
\[
f(\alpha^R)\to+\infty
\qquad\text{as}\qquad
\alpha^R\to\pm\infty.
\]
At the endpoints of the interval $(1-M^R,M^R-\eta)$, one has
\[
f(1-M^R)=-1,
\qquad
f(M^R-\eta)=-1.
\]
Therefore $f$ has at least one real root in $(-\infty,1-M^R)$ and at least one real root in
$(M^R-\eta,\infty)$.

Now choose
\[
\alpha_*^R
:=
M^R\frac{1-\eta^{2/3}}{1+\eta^{2/3}}.
\]
By \eqref{nonlinear-stability1}, we have
\[
\alpha_*^R\in(1-M^R,M^R-\eta).
\]
Moreover,
\[
U_*=\alpha_*^R+M^R
=
\frac{2M^R}{1+\eta^{2/3}},
\qquad
V_*=\frac{\alpha_*^R-M^R}{\eta}
=
-\frac{2M^R}{\eta^{1/3}(1+\eta^{2/3})}.
\]
Using again \eqref{nonlinear-stability1}, one obtains
\[
U_*>\sqrt{1+\eta^{2/3}},
\qquad
|V_*|>\sqrt{1+\eta^{-2/3}},
\]
and hence
\[
(U_*^2-1)(V_*^2-1)>1.
\]
Since
\[
(U_*^2-1)(V_*^2-1)=f(\alpha_*^R)+1,
\]
it follows that
\[
f(\alpha_*^R)>0.
\]
Together with
\[
f(1-M^R)=-1,
\qquad
f(M^R-\eta)=-1,
\]
the intermediate value theorem implies that $f$ has at least two real roots in
\[
(1-M^R,M^R-\eta).
\]

Since $f$ is a quartic polynomial and already has one real root on each side of this interval,
it follows that $f$ has exactly four real roots in total. Denote them by
\[
\alpha_3^R<\alpha_4^R<\alpha_5^R<\alpha_6^R.
\]
Moreover, exactly two of them lie in $(1-M^R,M^R-\eta)$, and the ordering is
\[
\alpha_3^R
<
1-M^R
<
\alpha_4^R
<
\alpha_5^R
<
M^R-\eta
<
\alpha_6^R.
\]

It remains to record the explicit expressions of the two roots in
\((1-M^R,M^R-\eta)\). Multiplying the equation \(f(\alpha^R)=0\)
by \(\eta^2\), we obtain the depressed quartic
\[
(\alpha^R)^4+p(\alpha^R)^2+q\alpha^R+r=0,
\]
where
\[
p:=-(2(M^R)^2+\eta^2+1),\qquad
q:=2M^R(1-\eta^2),\qquad
r:=(M^R)^4-(M^R)^2(\eta^2+1).
\]
Let \(y_*\) be the largest real root of
\[
\bigl((p+y_*)^2-r\bigr)(p+2y_*)=\frac{q^2}{4}.
\]
Since the quartic has four distinct real roots, the Ferrari resolvent roots
correspond to the three pairings of these roots. Hence the largest resolvent
root satisfies \(p+2y_*>0\). Then Ferrari's formula gives four roots of the quartic. In view of the ordering already proved, the two roots in \((1-M^R,M^R-\eta)\) are
\begin{align}
\alpha_4^R
&=
\frac12\left(
-\sqrt{p+2y_*}
+
\sqrt{-\left(3p+2y_*-\frac{2q}{\sqrt{p+2y_*}}\right)}
\right),
\nonumber\\
\alpha_5^R
&=
\frac12\left(
\sqrt{p+2y_*}
-
\sqrt{-\left(3p+2y_*+\frac{2q}{\sqrt{p+2y_*}}\right)}
\right).
\nonumber
\end{align}

\end{proof}

\bibliography{aipsamp}

\end{document}